# A New Approximate Solution of Time-Fractional, Non-linear Schrödinger Equations Using Fractional Reduced Differential Transformation

**Brajesh Kumar Singh[1] and Pramod Kumar[2,*]**

[1,2]Department of Applied Mathematics, School for Physical sciences,
Babasaheb Bhiamrao Ambedkar University, Lucknow 226025 India
{[1]bksingh0584, [2]bbaupramod}@gmail.com

**Abstract.** This paper is concerned with an alternative analytical solution of time-fractional nonlinear Schrödinger equation and nonlinear coupled Schrödinger equation obtained by employing fractional reduced differential transform method. The proposed solutions are obtained in series form, converges to the exact solution very rapidly. These results are agreed well with the results obtained by using differential transform method, homotopy perturbation method, homotopy analysis method and Adomian decomposition method. However, the computations shows that the described method is easy to apply, and it needs small size of computation contrary to the existing above said methods.

**Keywords:** Nonlinear (coupled) Schrödinger equations, fractional reduced differential transform, Mittag-Leffler function, analytic solution

**Mathematics Subject Classification 2010:** 26A33 · 34A08 · 34K37 · 34M25

## 1. Introduction

This paper focuses on an alternate analytical solution one dimensional (non)linear Schrödinger equation of the form:

$$i D_t^\alpha u + u_{xx} + \sigma |u|^2 u = 0, \qquad u(x,0) = f(x), \qquad x \in R, t \geq 0 \quad (1.1)$$

And one dimensional coupled nonlinear Schrödinger equation of the form:

$$\begin{cases} i D_t^\alpha u + u_{xx} + \sigma\left(|u|^2 + |v|^2\right)u = 0, & u(x,0) = f(x), \\ i D_t^\alpha v + v_{xx} + \sigma\left(|u|^2 + |v|^2\right)v = 0, & v(x,0) = g(x). \end{cases}, x \in R, t \geq 0 \quad (1.2)$$

Here $\sigma$ is a constant, $u(x,t)$ and $v(x,t)$ complex function, $i = \sqrt{-1}$, and $f, g$ are smooth functions. For $\sigma = 0$, Eq. (1.1) and Eq. (1.2) represent linear model of Schrödinger equation (LSE) and coupled Schrödinger equation. In the recent years, a great effort has been made in finding the exact solution of (non)linear differential equations to understand the most (non)linear physical phenomena. The fractional model of nonlinear Schrödinger equation (NLSE) is one of the most efficient universal models which describe various physical nonlinear systems. For example, NLSE is used to describe the evolution of slowly varying packets of quasi monochromatic waves in weakly nonlinear media with dispersion. Unlike the LSE, the NLSE never describes the time evolution of a quantum state. NLSE has found its various applications in nonlinear wave propagation in dispersive and inhomogeneous media such as: dynamics in particle accelerators [1], non-uniform dielectric media, solitary waves in piezoelectric semiconductors, hydrodynamics and plasma waves, nonlinear optical waves, quantum condensates, nonlinear acoustics, heat pulses in solids and various other nonlinear instability phenomena, and mean field theory of Bose–Einstein condensates [2–5].

The fractional model of Schrödinger equations has been solved by various methods [7, 14], among others, homotopy perturbation method [6,8,18], Adomian decomposition method [8,11], two dimensional differential transform methods [10], differential transform method [16], variational iteration method [17], fractional Riccati expansion method [15],Wang and Xu [9] used integral transforms technique to solve the space time fractional Schrodinger equation. Wang [12] employed split-step finite difference method to solve the nonlinear Schrödinger equations. An extensive work has been carried out by Masemola et al. [13] who proposed optical solitons and conservation laws for driven nonlinear Schrödinger's equation with linear attenuation and detuning. Recently, the fractional model of (coupled) nonlinear Schrödinger's equation has been solved by using HAM: homotopy analysis method [19,23] by Bakkyaraj and Sahadevan, Jacobi spectral collocation method by Bhrawy et al. [20], linearly implicit conservative difference scheme by Wang et al. [21] and Kudryashov method by Eslami [22].

This rest of the paper is sketched as follows: the basic definitions and notations of fractional calculus theory are revisited in Section 2 and that of fractional reduced differential transform method (FRDTM) in Section 3. The main goal of this paper is to present analytical solution of one dimensional (coupled) linear and nonlinear Schrödinger equation in Section 4. Finally, the concluding remark to the proposed results is given in Section 5.

2. **Fractional Calculus Theory**

This section deals with the basic of fractional calculus based on Liouville [25], which is needed to complete our study.

**Definition 2.1** Let $\mu \in \mathbb{R}$ and $m \in \mathbb{N}$. Then $f: \mathbb{R}^+ \to \mathbb{R}$ belongs to the space $C_\mu$ if there exists $k \in \mathbb{R}, k > \mu$ and $g \in C[0,\infty)$ such that $f(x) = x^q g(x)$. Moreover, $f \in C_\mu^m$ if $f^{(m)} \in C_\mu$.

**Definition 2.2** Let $J_x^\alpha$ $(\alpha \geq 0)$ be the Riemann Liouville fractional integral operator [25], then

$$J_x^\alpha f(x) = \frac{1}{\Gamma(\alpha)} \int_0^x (x-t)^{\alpha-1} f(t) dt, \ \alpha > 0, x > 0 \ \ and \ \ J_x^0 f(x) = f(x).$$

(2.1)

Where $\Gamma(z) := \int_0^\infty e^{-t} t^{z-1} dt, z \in \mathbb{C}$.

Let $m-1 < \alpha \leq m$ and $m \in \mathbb{N}$, $f \in C_\mu^m$ ($\mu \geq -1$), $\alpha, \beta \geq 0$ and $\gamma > -1$, then the operator $J_x^\alpha$ satisfy the following properties:

i) $J_x^\alpha J_x^\beta f(x) = J_x^{\alpha+\beta} f(x) = J_x^\beta J_x^\alpha f(x)$,

ii) $J_x^\alpha x^\gamma = \dfrac{\Gamma(1+\gamma)}{\Gamma(1+\gamma+\alpha)} x^{\gamma+\alpha}$, x>0.

Moreover, Caputo and Mainardi [26] developed a modified fractional differentiation operator $D_x^\alpha$ to overcome the discrepancy of Riemann-Liouville derivative.

**Definition 2.3** Let $m-1 < \alpha \leq m$, $m \in \mathbb{N}, x > 0$, and then Caputo fractional derivative of $f \in C_\mu$ [26] is read as

$$D_x^\alpha f(x) = J_x^{m-\alpha} D_x^m f(x) = \frac{1}{\Gamma(m-\alpha)} \int_0^x (x-t)^{m-\alpha-1} f^{(m)}(t) dt. \quad (2.2)$$

The basic properties of $D_x^\alpha$ are as follows:

**Lemma 2.1** If $m-1 < \alpha \leq m$, $m \in N$ and $f \in C_\mu^m$, $\mu \geq -1$, then

a) $D_x^\alpha D_x^\beta f(x) = D_x^{\alpha+\beta} f(x) = D_x^\beta D_x^\alpha f(x)$,

b) $D_x^\alpha x^\gamma = \dfrac{\Gamma(1+\gamma)}{\Gamma(1+\gamma-\alpha)} x^{\gamma-\alpha}$, x>0

c) $D_x^\alpha J_x^\alpha f(x) = f(x)$, x>0,

d) $J_x^\alpha D_x^\alpha f(x) = f(x) - \sum_{k=0}^m f^{(k)}(0^+) \dfrac{x^k}{k!}$, x>0,

For details study of fractional derivatives we refer the readers to [25-27].

## 3. FRDT method

This section concerned with the discussion of some basic results as in [28-32], on fractional reduced differential transform to complete the paper. Throughout the paper, we denote the original function by $\phi(x,t)$ (lowercase) while it's fractional reduced differential transform (FRDT) by $\Phi_k(x,t)$ (uppercase).

**Definition 3.1** FRDT (spectrum) of an analytic and continuously differentiable function $w(x,t)$ is defined by

$$W_k(x) = \frac{1}{\Gamma(k\alpha+1)} \{D_t^{\alpha k} w(x,t)\}_{t=t_0} \quad (3.1)$$

where $\alpha$ is order of fractional derivative. The inverse FRDT of $W_k(x)$ is defined as follows

$$w(x,t) = \sum_{k=0}^{\infty} W_k(x)(t-t_0)^{k\alpha}. \quad (3.2)$$

From Eq. (3.2) and (3.3), one get

$$w(x,t) = \sum_{k=0}^{\infty} \frac{1}{\Gamma(k\alpha+1)} \{D_t^{\alpha k} w(x,t)\}_{t=t_0} (t-t_0)^{k\alpha}. \quad (3.3)$$

In particular for $t_0 = 0$, we get

$$w(x,t) = \sum_{k=0}^{\infty} W_k(x) t^{k\alpha} = \sum_{k=0}^{\infty} \frac{1}{\Gamma(k\alpha+1)} \{D_t^{\alpha k} w(x,t)\}_{t=0} t^{k\alpha}. \quad (3.4)$$

**Definition 3.2** The Mittag-Leffler function $E_\alpha(z)$ with $\alpha > 0$ is defined by

$$E_\alpha(z) = \sum_{k=0}^{\infty} \frac{z^k}{\Gamma(1+k\alpha)}, \quad (3.5)$$

is valid in the whole complex plane, and is an advanced form of $\exp(z)$. $\exp(z) = \lim_{\alpha \to 1} E_\alpha(z)$.

**Theorem 3.1** Let $U_k(x)$ and $V_k(x)$ be spectrum of analytic and continuously differentiable function $u(x,t)$ and $v(x,t)$ respectively, then

a) If $w(x,t) = u(x,t)v(x,t)$, then
$$W_k(x) = U_k(x) \otimes V_k(x) = \sum_{r=0}^{k} U_r(x) V_{k-r}(x).$$

b) If $w(x,t) = \ell_1 u(x,t) \pm \ell_2 v(x,t)$, then
$$W_k(x) = \ell_1 U_k(x) \pm \ell_2 V_k(x).$$

c) If $\psi(x,t) = u(x,t)v(x,t)w(x,t)$, then

$$\Psi_k(x) = U_k(x) \otimes V_k(x) \otimes W_k(x) = \sum_{r=0}^{k} \sum_{i=0}^{r} U_i(x) \, V_{r-i}(x) W_{k-r}(x).$$

d) If $\psi(x,t) = D_t^{N\alpha} u(x,t)$, then $\Psi_k(x) = \dfrac{\Gamma(1+(k+N)\alpha)}{\Gamma(1+k\alpha)} U_{k+N}(x).$

e) If $\theta(x,t) = x^m t^n \psi(x,t)$, then

$$\Theta_k(x) = \begin{cases} x^m \Psi_{k\alpha-n}(x), & \text{if } k\alpha \geq n \\ 0, & \text{else.} \end{cases}$$

f) If $\theta(x,t) = x^m t^n$, then $\Theta_k(x) = x^m \delta(k\alpha - n),$ where the function $\delta$ is defined by

$$\delta(k) = \begin{cases} 1 & \text{if } k = 0 \\ 0 & \text{otherwise} \end{cases}.$$

## 4. Numerical Examples

The main goal of the paper is to present the analytical computation of five test problems of fractional model of (N)LSE and coupled NLSE.

**Problem 4.1** Consider one dimensional LSE (1.1) with $\sigma = 0$, $f(x) = 1 + \cosh ax$, $a$ is constant. The FRDT of LSE (1.1) with $\sigma = 0$ produces the following recurrence relation:

$$\frac{\Gamma(1+(1+k))\alpha}{\Gamma(1+k\alpha)} U_{k+1}(x) = i \frac{\partial^2 U_k(x)}{\partial x^2}, \quad U_0(x) = 1 + \cosh ax, \ k \geq 1.$$

On solving the above recurrence relation, we get

$$U_1(x) = \frac{1}{\Gamma(1+\alpha)} i a^2 \cosh ax;$$

$$U_2(x) = \frac{1}{\Gamma(1+2\alpha)} i^2 a^4 \cosh ax;$$

$$U_3(x) = \frac{1}{\Gamma(1+3\alpha)} i^3 a^6 \cosh ax;$$

$$U_4(x) = \frac{1}{\Gamma(1+4\alpha)} i^4 a^8 \cosh ax;$$

$$U_5(x) = \frac{1}{\Gamma(1+5\alpha)} i^5 a^{10} \cosh ax;$$

$$\vdots \qquad \vdots$$

(4.1)

The inverse FRDT (3.5) leads to

$$u(x,t) = \sum_{k=0}^{\infty} U_k(x) t^{k\alpha} = U_0(x) + U_1(x) t^\alpha + U_2(x) t^{2\alpha} + U_3(x) t^{3\alpha} + \ldots$$

$$= 1 + \cosh ax \left( 1 + \frac{1}{\Gamma(1+\alpha)} i a^2 t^\alpha + \frac{1}{\Gamma(1+2\alpha)} i^2 a^4 t^{2\alpha} + \cdots + \frac{1}{\Gamma(1+2k\alpha)} i^k a^{2k} t^{k\alpha} + \ldots \right)$$

(4.2)

Consequently, the series in the closedform

$$u(x,t) = 1 + \cosh ax \sum_{k=0}^{\infty} \frac{(ia^2 t^\alpha)^k}{\Gamma(1+k\alpha)}. \quad (4.3)$$

In particular, the approximate solution (4.3) for $\alpha = 1$, $a = 2$ is exactly same as obtained by using ADM [8], DTM [10], HPM [8, 19] and HAM [23]. The physical behavior of real part $\phi(x,t)$ and imaginary part $\psi(x,t)$ of $u(x,t)$ with $\alpha = 0.5, 0.9$ for different domain is depicted in Fig 4.1.

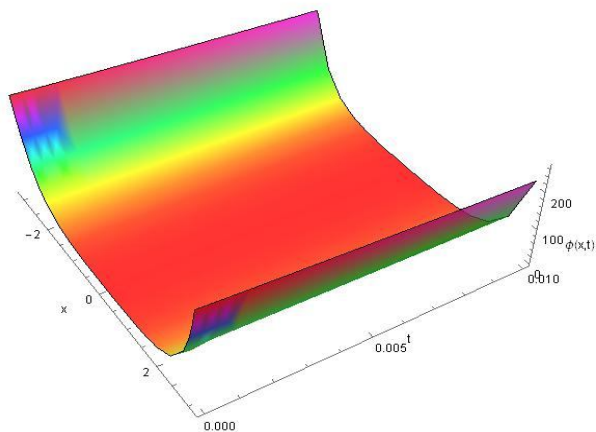
$\phi$ : for $\alpha = 0.9$ in the domain $x \in (-\pi, \pi), t \in (0, 0.01)$

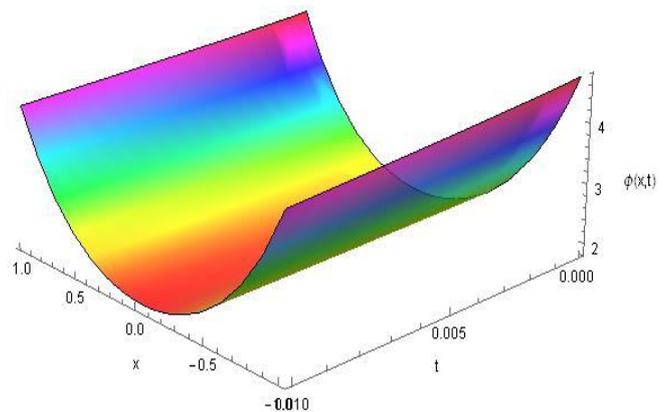
$\phi$ : for $\alpha = 0.9$ in the domain $t \in (0, 0.01)$, $x \in (-1,1)$

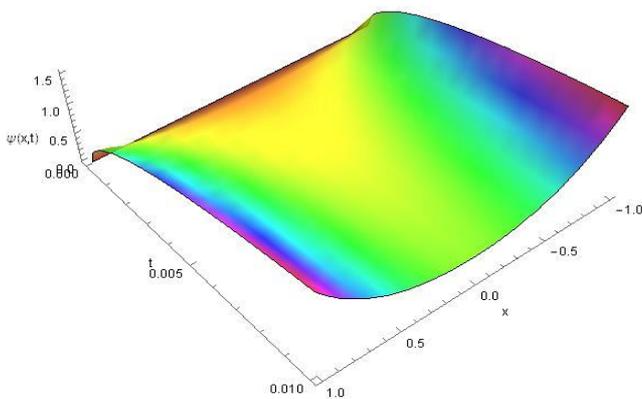
$\psi$ : for $\alpha = 0.5$ in the domain $t \in (0, 0.01)$ and $x \in (-1,1)$

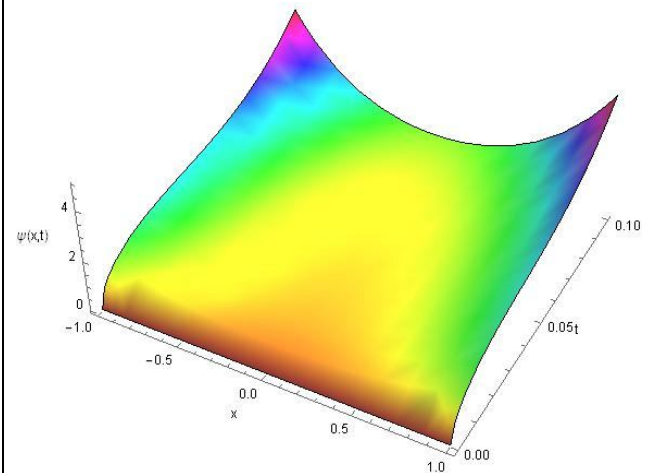
$\psi$ : for $\alpha = 0.5$ in the domain $t \in (0,1)$ and $x \in (-1,1)$

**Fig. 4.1** Phisical solution behavior of $\phi$ and $\psi$ at different time levels in different domain

**Problem 4.2** Consider one dimensional LSE (1.1) with $\sigma = 0$, $f(x) = e^{inx}$, $n$ is constant

The FRDT of LSE (1.1) with $\sigma = 0$ produces the following recurrence relation:

$$\frac{\Gamma(1+(1+k))\alpha}{\Gamma(1+k\alpha)}U_{k+1}(x)=i\frac{\partial^2 U_k(x)}{\partial x^2}, \quad U_0(x)=e^{inx}, k\geq 1.$$

On solving the above recurrence relation, we get

$$U_1(x)=\frac{-n^2 i e^{nix}}{\Gamma(1+\alpha)}; U_2(x)=-\frac{n^4 e^{nix}}{\Gamma(1+2\alpha)}; U_3(x)=\frac{n^6 i e^{nix}}{\Gamma(1+3\alpha)}; U_4(x)=\frac{n^8 e^{nix}}{\Gamma(1+4\alpha)};\cdots$$

The inverse FRDT method (3.5) leads to

$$u(x,t)=e^{nix}-i\frac{n^2}{\Gamma(1+\alpha)}e^{nix}t^\alpha-\frac{n^4}{\Gamma(1+2\alpha)}e^{nix}t^{2\alpha}+i\frac{n^6}{\Gamma(1+3\alpha)}e^{nix}+\frac{n^8}{\Gamma(1+4\alpha)}e^{nix}+\cdots$$
$$=e^{nix}\sum_{k=0}^{\infty}\frac{\left(n^2(-i)t^\alpha\right)^k}{\Gamma(1+k\alpha)}.$$

(4.4)

The physical behavior of real and imaginary parts of $u(x,t)$ is depicted in Fig 4.2.

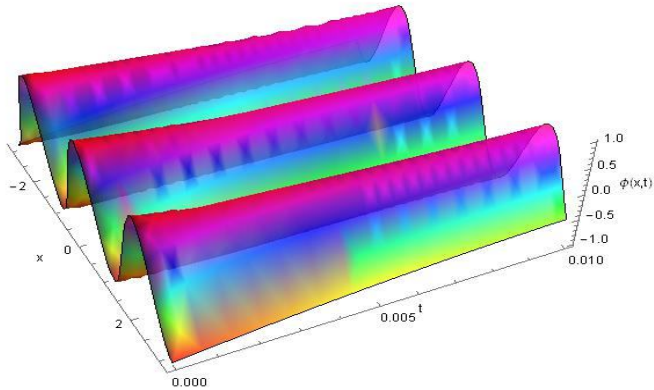
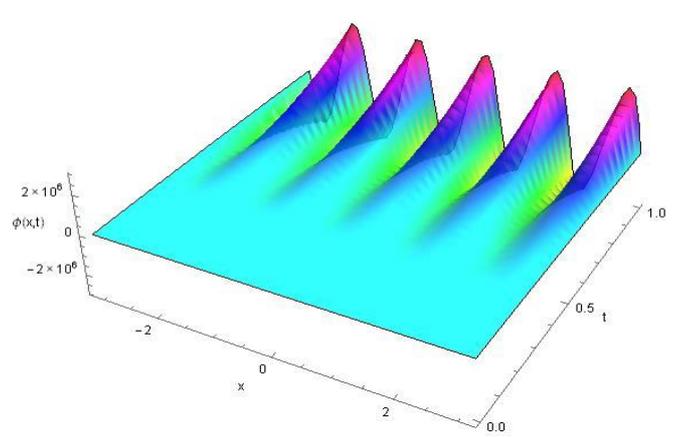

$\phi$: for $\alpha=0.5, n=3$ in the domain $x\in(-\pi,\pi), t\in(0,0.01)$ | $\phi$: for $\alpha=0.5, n=5$ in the domain $x\in(-\pi,\pi), t\in(0,1)$

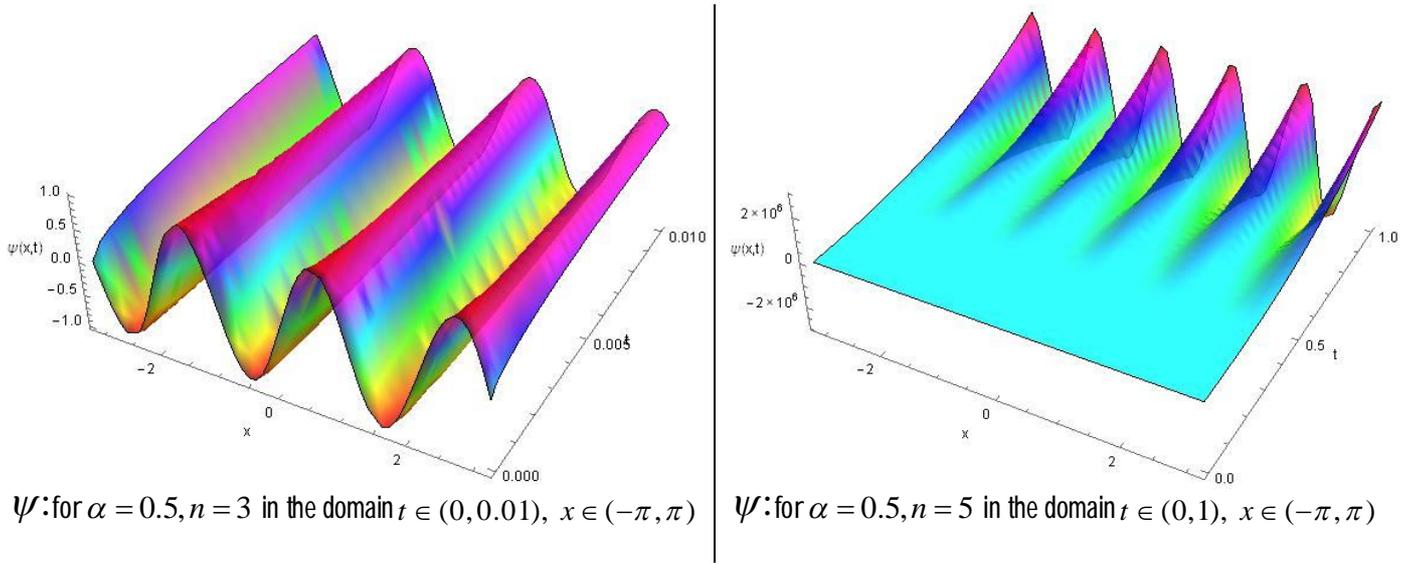

$\psi$: for $\alpha = 0.5, n = 3$ in the domain $t \in (0, 0.01),\ x \in (-\pi, \pi)$ | $\psi$: for $\alpha = 0.5, n = 5$ in the domain $t \in (0,1),\ x \in (-\pi, \pi)$

**Fig. 4.2** Physical solution behavior of $\phi$ and $\psi$ at different time levels in different domain

**Problem 4.3** Consider one dimensional NLSE (1.1) with $f(x) = e^{inx}$, $n$ is constant

In this case, FRDT of LSE (1.1) produces the following recurrence relation:

$$\frac{\Gamma(1+(1+k))\alpha}{\Gamma(1+k\alpha)} U_{k+1}(x) = i\left(\frac{\partial^2 U_k(x)}{\partial x^2} + \sigma \sum_{k_2=0}^{k}\sum_{k_1=0}^{k_2} \bar{U}_{k_1}(x)\, U_{k_2-k_1}(x) U_{k-k_2}(x)\right),\ k \geq 1$$

$$U_0(x) = e^{inx}$$

On solving the above relation, we get

$$U_1(x) = i\, \frac{1}{\Gamma(1+\alpha)}(\sigma - n^2)\, e^{inx}$$

$$U_2(x) = i^2\, \frac{1}{\Gamma(1+2\alpha)}(\sigma - n^2)^2\, e^{inx}$$

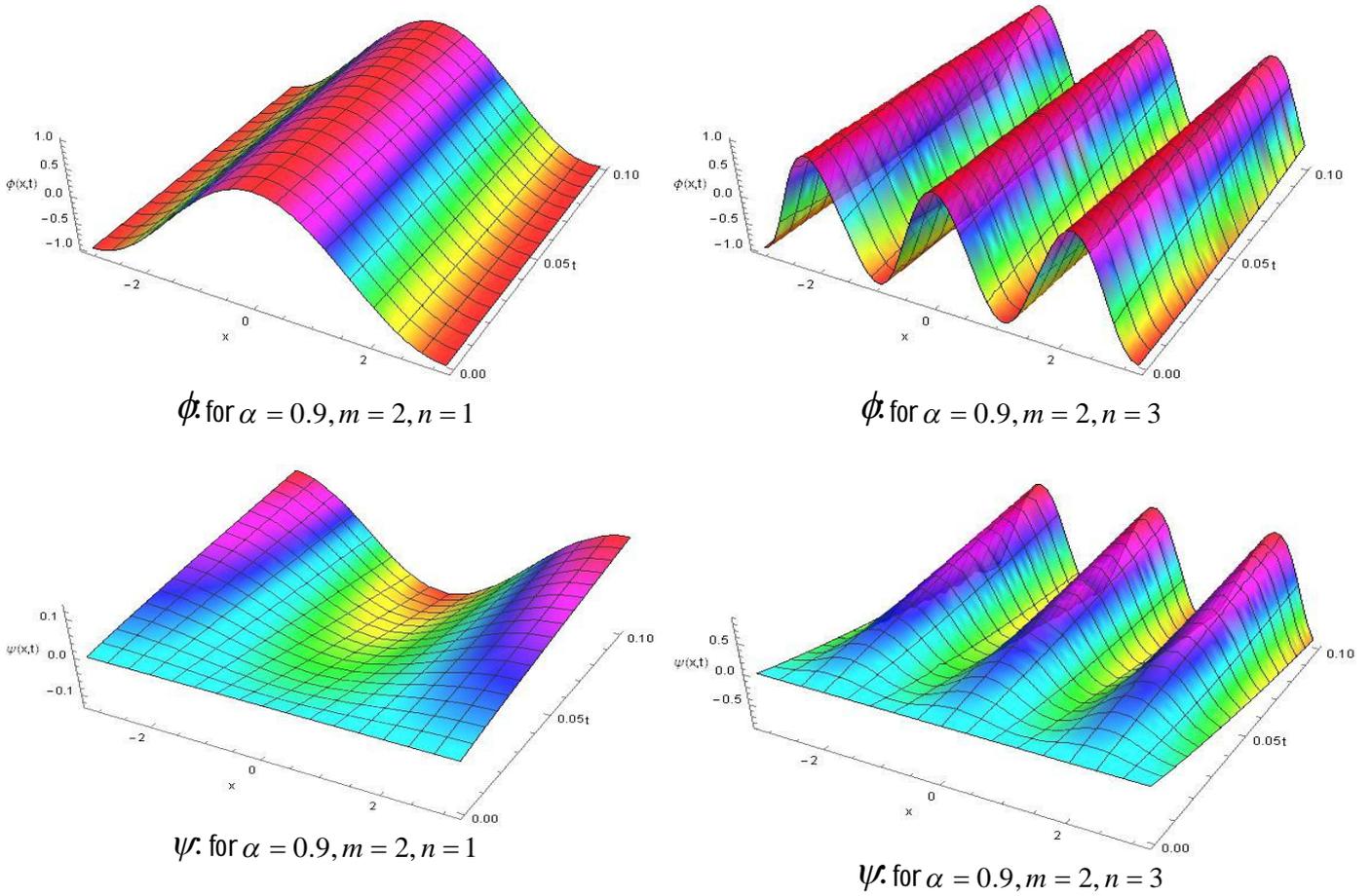

Fig. 4.3 Physical solution behavior of $\phi$ and $\psi$ at different time levels in the domain $x \in (-\pi, \pi), t \in (0, 0.1)$

$$U_3(x) = -i^3 \frac{(\sigma - n^2)^2}{\Gamma(1+3\alpha)} \left[ (n^2 - 3\sigma) + \frac{\sigma\Gamma(1+2\alpha)}{\Gamma(1+3\alpha)} \right] e^{inx}$$

$$U_4(x) = i^4 \frac{(-\sigma + n^2)^3}{\Gamma(1+4\alpha)} \left[ (-3\sigma + n^2) + \frac{\sigma\Gamma(1+2\alpha)}{\Gamma(1+\alpha)^2} - \frac{2\sigma\Gamma(1+3\alpha)}{\Gamma(1+\alpha)\Gamma(1+2\alpha)} + \frac{\sigma\Gamma(1+3\alpha)}{\Gamma(1+\alpha)^3} \right] e^{inx}$$

$\vdots$

The inverse FRDT leads to

$$u(x,t) = e^{inx} + i\frac{1}{\Gamma(1+\alpha)}(\sigma - n^2) e^{inx} t^\alpha + i^2 \frac{1}{\Gamma(1+2\alpha)}(\sigma - n^2)^2 e^{inx} t^{2\alpha} + \ldots$$

For $\alpha = 1$, the above series solution reduces to the closed form

$$u(x,t) = \sum_{k=0}^{\infty} \frac{\left(i(m-n^2)t\right)^k}{\Gamma(1+k)} e^{inx} = e^{i(nx+(\sigma-n^2)t)}. \tag{4.5}$$

Moreover, the approximate solution (4.5) with $\alpha = 1, \sigma = -2, n = 1$ is exactly same as obtained by using DTM [10], HPM [8], VIM [17]. The physical behavior of real part $\phi(x,t)$ and imaginary part $\psi(x,t)$ of $u(x,t)$ with $\alpha = 0.9, \sigma = 2, n = 1,3$ in the domain $t \in (0,0.1)$ and $x \in (-\pi, \pi)$ is depicted in Fig 4.3.

**Problem 4.4** Now, consider the fractional model of NLSE with trapping potential

$$i D_t^\alpha u = -\frac{1}{2} u_{xx} + u \cos^2 x + |u|^2 u, \tag{4.6}$$
$$u(x,0) = \sin x$$

The FRDT method, Eq. (5.6) reduces to a set of recurrence relation as follows:

$$\frac{\Gamma(1+(1+k))\alpha}{\Gamma(1+k\alpha)} U_{k+1}(x) = i \left[ \frac{1}{2} \frac{\partial^2 U_k}{\partial x^2} - U_k \cos^2 x - \sum_{k_2=0}^{k} \sum_{k_1=0}^{k_2} \bar{U}_{k_1}(x) U_{k_2-k_1}(x) U_{k-k_2}(x) \right] \tag{4.7}$$
$$U_0 = \sin x$$

On solving the recurrence relation (4.7), we get

$$U_0 = \sin x, \ U_1 = -\frac{3i \sin x}{2\Gamma(1+\alpha)}, \ U_2 = \frac{9i^2}{4} \frac{\sin x}{\Gamma(1+2\alpha)}, \ U_3 = \frac{9i^3}{8} \frac{\sin x}{\Gamma(1+3\alpha)} \left[ (-5 + 2\cos 2x) + \frac{2\Gamma(1+2\alpha) \sin^2 x}{\Gamma(1+\alpha)^2} \right]$$

$$U_4 = -\frac{9i^4 \sin x}{16 \Gamma(1+4\alpha)} \left[ -7 + 22\cos 2x - \frac{(1+11\cos 2x)\Gamma(1+2\alpha)}{\Gamma(1+\alpha)^2} - \frac{12\Gamma(1+3\alpha) \sin^2 x}{\Gamma(1+\alpha)\Gamma(1+2\alpha)} + \frac{6\Gamma(1+3\alpha) \sin^2 x}{\Gamma(1+\alpha)^3} \right]$$
$$\vdots \qquad \vdots$$

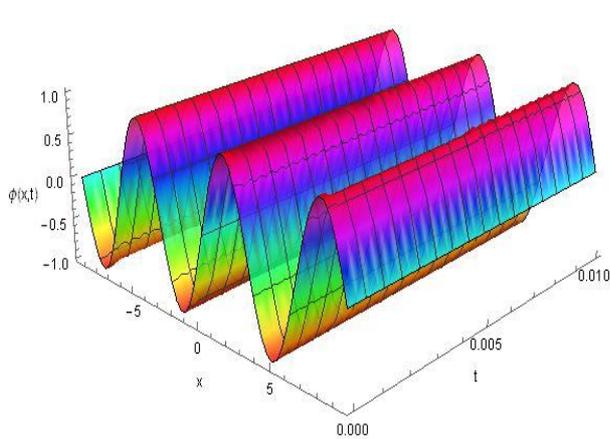

$\phi: x \in (-3\pi, 3\pi)$

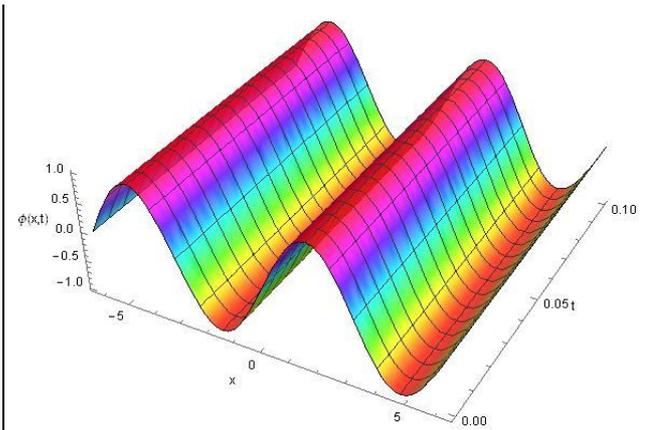

$\phi: x \in (-2\pi, 2\pi)$

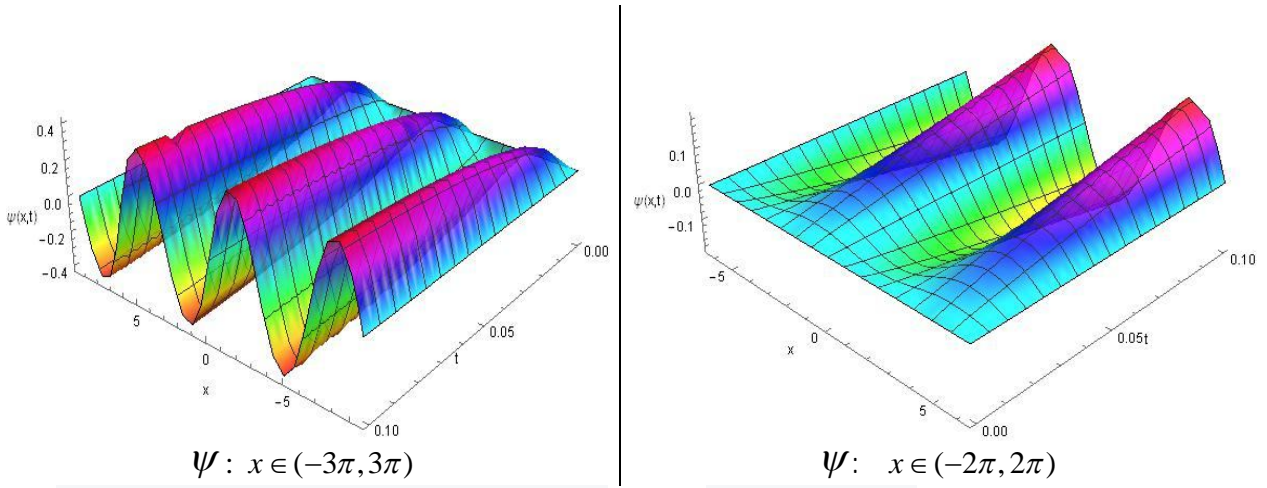

$\psi:\ x\in(-3\pi,3\pi)$  |  $\psi:\ x\in(-2\pi,2\pi)$

**Fig. 4.4** Physical behavior of $\phi$ and $\psi$ for $\alpha=0.5$ at in the domain $x\in(-\pi,\pi), t\in(0,0.1)$

The inverse FRDT leads to

$$u(x,t)=\sin x-\frac{3i}{2}\frac{\sin x}{\Gamma(1+\alpha)}t^{\alpha}+\frac{9i^2}{4}\frac{\sin x}{\Gamma(1+2\alpha)}t^{2\alpha}$$
$$+\frac{9i^3}{8}\frac{\sin x}{\Gamma(1+3\alpha)}\left[-5+2\cos 2x+\frac{2\Gamma(1+2\alpha)\sin^2 x}{\Gamma(1+\alpha)^2}\right]t^{3\alpha}+\ldots \quad (4.8)$$

The approximate solution (4.8) with is $\alpha=1$ reduces to

$$u(x,t)=e^{\left(-\frac{3it}{2}\right)}\sin x. \quad (4.9)$$

This is exactly same as obtained by DTM [10]. Fig 4.4 depicts the physical behavior of $\phi(x,t)$ and $\psi(x,t)$ for $\alpha=0.5$ at $t\in(0,0.1)$ in the domain $x\in(-2\pi,2\pi)$ or $x\in(-3\pi,3\pi)$.

**Problem 4.5** Consider fractional model of coupled NLSE (1.2) with $\sigma=2$, $f(x)=ae^{inx}$ and $g(x)=be^{imx}$ as in [19]:

FRDT to Eq. (1.2) with $\sigma=2$, $f(x)=ae^{inx}$ and $g(x)=be^{imx}$ produces the following set of recurrence relation:

$$\frac{\Gamma(1+(1+k))\alpha}{\Gamma(1+k\alpha)}U_{k+1}(x) = \frac{i}{2}\frac{\partial^2 U_k}{\partial x^2}$$

$$+2i\left(\sum_{k_2=0}^{k}\sum_{k_1=0}^{k_2}\overline{U}_{k_1}(x)U_{k_2-k_1}(x)U_{k-k_2}(x) + \sum_{k_2=0}^{k}\sum_{k_1=0}^{k_2}\overline{V}_{k_1}(x)V_{k_2-k_1}(x)U_{k-k_2}(x)\right), U_0 = ae^{inx}$$

$$\frac{\Gamma(1+(1+k))\alpha}{\Gamma(1+k\alpha)}V_{k+1}(x) = \frac{i}{2}\frac{\partial^2 V_k}{\partial x^2}$$

$$+2i\left(\sum_{k_2=0}^{k}\sum_{k_1=0}^{k_2}\overline{U}_{k_1}(x)U_{k_2-k_1}(x)V_{k-k_2}(x) + \sum_{k_2=0}^{k}\sum_{k_1=0}^{k_2}\overline{V}_{k_1}(x)V_{k_2-k_1}(x)V_{k-k_2}(x)\right), V_0 = be^{imx}.$$

(4.10)

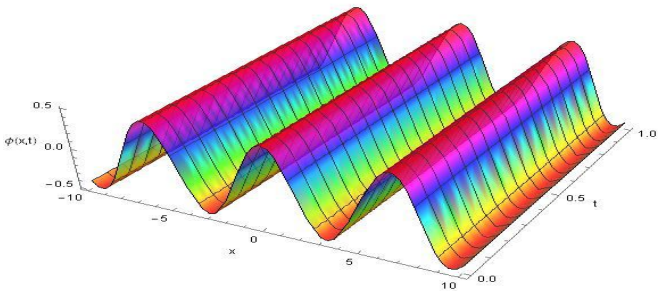

$\phi_1 : a = b = 0.5, n = 1, m = 1.5$

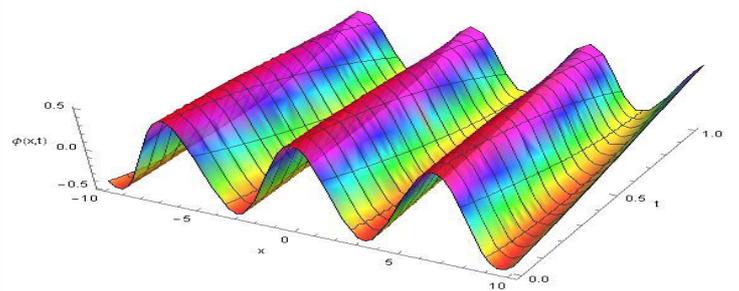

$\phi_1 : a = 0.5, n = 1 = b, m = 1.5$

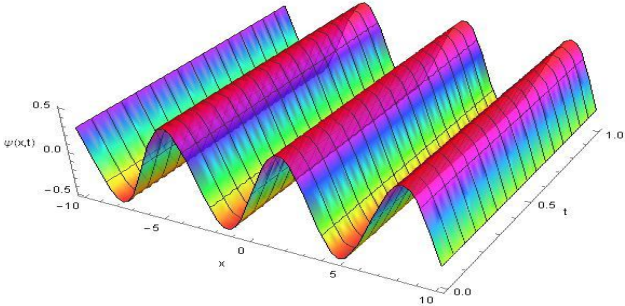

$\psi_1 : a = b = 0.5, n = 1, m = 1.5$

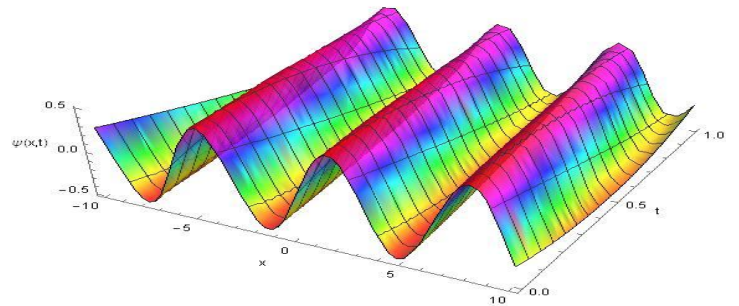

$\psi_1 : a = 0.5, b = n = 1, m = 1.5$

**Fig. 4.5** Physical behavior of real and imaginary parts $u = \phi_1 + \psi_1$ for $\alpha = 0.9$ at $t = (0,1)$ in the domain $x \in (-10,10)$

On solving the recurrence relation (4.10) for $k = 1, 2, \ldots$, we get

$$U_1 = \frac{iae^{inx}(2a^2 + 2b^2 - n^2)}{\Gamma(1+\alpha)}, \quad V_1 = \frac{ibe^{imx}(2a^2 + 2b^2 - m^2)}{\Gamma(1+\alpha)}$$

$$U_2 = -\frac{ae^{inx}(2a^2 + 2b^2 - n^2)^2}{\Gamma(1+2\alpha)}, V_2 = -\frac{be^{imx}(2a^2 + 2b^2 - m^2)^2}{\Gamma(1+2\alpha)}$$

⋮        ⋮

The inverse FRDT leads to

$$u(x,t) = ae^{inx} + \frac{iae^{inx}(2a^2+2b^2-n^2)t^\alpha}{\Gamma(1+\alpha)} - \frac{ae^{inx}(2a^2+2b^2-n^2)^2 t^{2\alpha}}{\Gamma(1+2\alpha)} + \cdots$$

$$v(x,t) = be^{imx} + \frac{ibe^{imx}(2a^2+2b^2-m^2)t^\alpha}{\Gamma(1+\alpha)} - \frac{be^{imx}(2a^2+2b^2-m^2)^2 t^{2\alpha}}{\Gamma(1+2\alpha)} + \cdots$$

(4.11)

Which In particular for $\alpha = 1,$ the approximate series solution (4.11) reduces to

$$u(x,t) = ae^{i\left(nx+(2a^2+2b^2-n^2)t\right)} \quad \text{and} \quad v(x,t) = be^{i\left(mx+(2a^2+2b^2-m^2)t\right)}.$$

The same solution is obtained by Tan and Boyd [24], Bakkyaraj and Sahadevan [19]. The physical behavior of real and imaginary parts of the solution $u = \phi_1 + \psi_1$ and $v = \phi_2 + \psi_2$ of (4.11) is depicted for $\alpha = 0.9,\ a = 0.5,\ n = 1, b = 0.5, 1,\ m = 1.5$ in $x \in (-10,10)$ at different time levels $t = (0,1)$

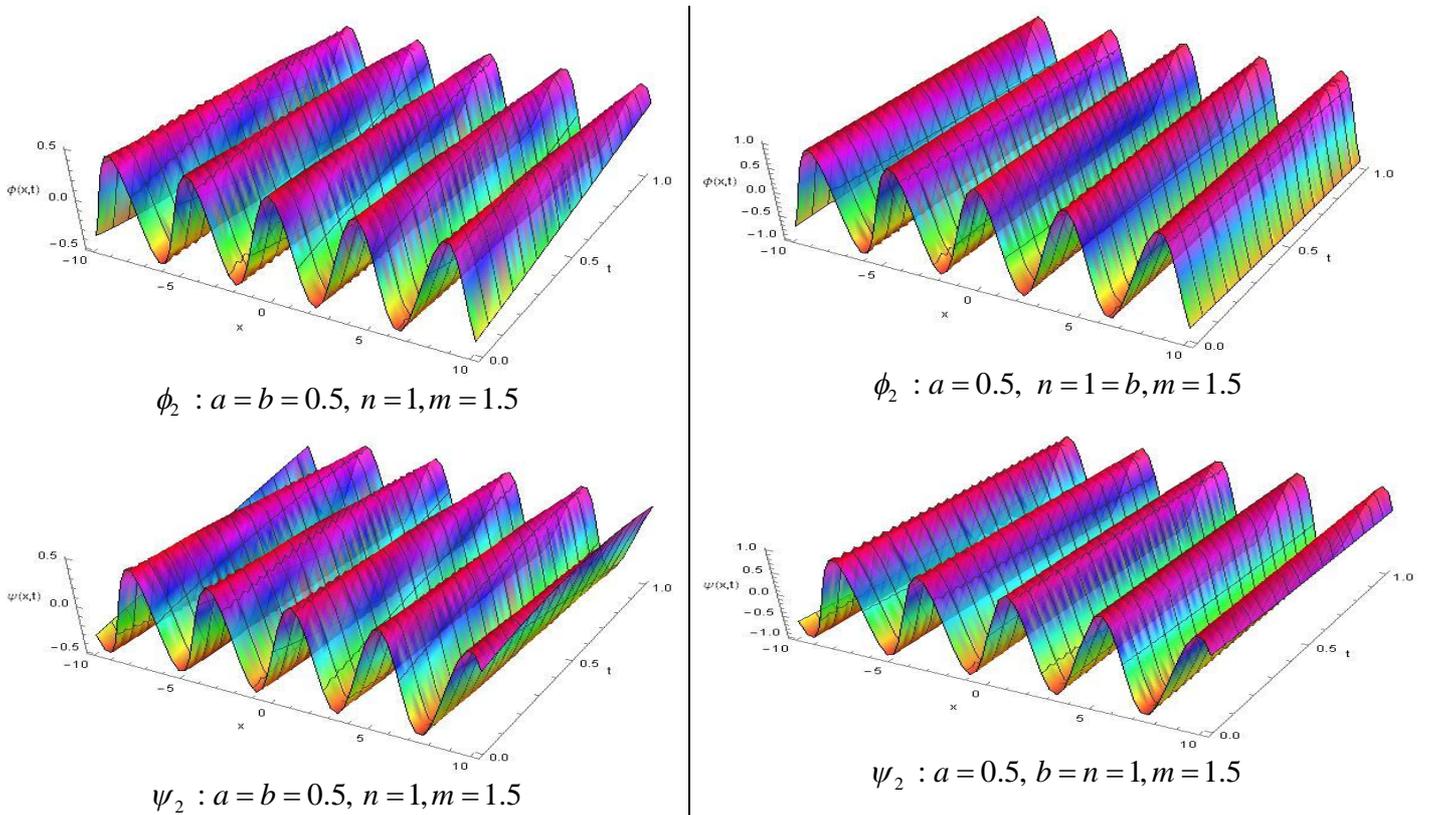

Fig. 4.6 Physical behavior of real and imaginary parts $v = \phi_2 + \psi_2$ for $\alpha = 0.9$ at $t = (0,1)$ in the domain $x \in (-10,10)$

## 5. Conclusion

In this paper, we illustrate the use of FRDTM in computation of an alternative analytical solution of time-fractional (non)linear Schrödinger equation, and coupled time-fractional

(non)linear Schrödinger equation. The FRDT solutions are obtained in infinite power series for appropriate initial condition, which converges to the exact solution rapidly. Moreover, the solutions are obtained without any perturbation, discretization or any other restrictive conditions. Five test problems are carried out to study the accurateness and effectiveness of the technique. The computed solutions are agreed well with those obtained by *Laplace transform method, differential transform method, homotopy perturbation method, homotopy analysis method* and *Adomain decomposition method*. However, the computations shows that the described method is easy to apply, and it needs small size of computation contrary to the existing above said methods.